\documentclass[a4paper,12pt]{article}
\usepackage{multicol,xcolor,hyperref}
\usepackage[english]{babel}
\usepackage{amsthm,amsmath,amssymb}
\usepackage{hyperref}

\newtheorem{rema}{Remark}
\newtheorem{lemma}{Lemma}

\newtheorem{conj}{Conjecture}

\newtheorem{thm}{Theorem}
\newtheorem{defn}{Definition}

\DeclareMathOperator{\card}{card}
\def\m{\diamond}
\def \a   {\alpha}
\def \b   {\beta}
\def \g   {\gamma}

\def \eps {\varepsilon}

\renewcommand{\i}{{\mathsf{i}}}

\def\E{{\mathbb{E}}}

\def\P{{\mathbb{P}}}

\def\F{{\cal{F}}}

\def\|{\,|\,}

\newcommand{\blue}[1]{\textcolor{blue}{#1}}

\title {Rigorous upper bound for the discrete Bak-Sneppen model}
\author{Stanislav Volkov\footnote{Centre for Mathematical Sciences, Lund University, Box 118 SE-22100, Lund, Sweden}} 
\begin {document}
\maketitle
\begin{abstract}
Fix some $p\in[0,1]$ and a positive integer $n$. The discrete Bak-Sneppen model is a Markov chain on the space of zero-one sequences of length $n$ with periodic boundary conditions. At each moment of time a minimum element (typically, zero) is chosen with equal probability, and it is then replaced alongside both its neighbours by independent Bernoulli($p$) random variables. Let $\nu^{(n)}(p)$ be the probability that an element of this sequence equals one under the stationary distribution of this Markov chain. It was shown in~\cite{BK} that $\nu^{(n)}(p)\to 1$ as $n\to\infty$ when $p>0.54\dots$; the proof there is, alas, not rigorous. The complimentary fact that $\displaystyle \limsup_{n\to\infty} \nu^{(n)}(p)< 1$ for $p\in(0,p')$ for some $p'>0$ is much harder; this was eventually shown in~\cite{MZ}.

The purpose of this note is to provide a {\em rigorous} proof of the result from~\cite{BK}, as well as to improve it, by showing that
$\nu^{(n)}(p)\to 1$  when $p>0.45$. (Our method, in fact, shows that with some finer tuning the same is true for $p>0.419533$.)
\end{abstract}



\section{Introduction}\label{Intro}
The classical Bak-Sneppen model~\cite{BAK,BAKS} is defined as a collection of $n$ individual species located equidistantly on a circumference, each possessing {\em a fitness}, which is a number in $(0,1)$. The process evolves in discrete time as follows. First, one finds  the node(s) with the minimal fitness (if there are more than one, each of them is chosen with equal probability), and then this individual is replaced by a new one, with a fitness drawn from a uniform $U(0,1)$ distribution. In such formulation, there are no interactions in the model,  and it is easy to see that the second highest fitness is always non-decreasing. Consequently, with a little extra work one can show that all but one  fitnesses converge to $1$ a.s. To make the model interesting, it is also assumed that the ``worst'' species is replaced together with both its immediate neighbours on the circumference, and each of the three new fitnesses is drawn independently from the same uniform distribution; as a result, the model becomes highly non-trivial. In particular, simulations indicate that  as time goes to infinity, for very large $n$ the distribution of each fitness converges to a uniform distribution with parameters $[f_c,1]$ where $f_c\approx 0.66$. To the best of our knowledge, this has not yet been shown rigorously.

The discrete version of the Bak-Sneppen model, proposed in~\cite{BK}, is defined as follows. Fix a positive integer $n\ge 3$ and $p\in[0,1]$ and consider a Markov chain $\xi(t)$, $t=0,1,2,\dots$, on the state space $\{0,1\}^n$ with the following transition probabilities. Let $\xi(t)=(x_0(t),x_1(t),\dots,x_{n-1}(t))$, $x_i(t)\in \{0,1\}$, and assume that $x_i(t)$ are the values assigned at time $t$ to $n$ vertices placed equidistantly on some circumference. Pick uniformly at random a vertex with index $i=i_t$ such that $x_i(t)=0$, and replace it\footnote{if all $x_i(t)=1$, then pick $i$ uniformly amongst $\{0,1,\dots,n-1\}$} and both of its neighbours (i.e., each of $\xi_{i-1}(t)$, $\xi_{i}(t)$, and $\xi_{i+1}(t)$) by an independent Bernoulli$(p)$ random variable, keeping all the remaining $\xi_j$ intact. Throughout the paper we assume periodic boundary conditions, that is, $n-1$ and $1$ are the neighbours of~$0$, and~$n-2$ and~$0$ are the neighbours of $n-1$, this is equivalent to addition/subtraction modulus~$n$. 

Formally, let $\bar\zeta(t)=\left(\zeta_{-1}(t),\zeta_{0}(t),\zeta_{1}(t),\hat\zeta(t)\right)$, $t=1,2,\dots$, be a collection of i.i.d.\ random vectors where $\zeta_k(t)$, $k=0,\pm1$, are Bernoulli$(p)$  and $\hat\zeta(t)$ are Uniform$[0,1]$, and the elements of $\bar\zeta(t)$ are also independent between themselves.  Assume $\xi(0)\in \{0,1\}^n$. The values of $\xi(t+1)$, $t=0,1,2,\dots$, are defined recursively as follows.

For $t\ge 1$, let $X_{t}$ is the number of zeroes in $\xi(t)$ and $i_{t}$ be the index of a randomly chosen zero\footnote{or a randomly chosen $1$, in case $\xi(t)$ does not contain any zeros} in the configuration $\xi(t)$. The locations of zeros in $\xi(t)$ are denoted by
$$
I(t):=\{\i_0(t),\i_1(t),\dots, \i_{X_{t}-1}(t)\}=
\{i:\ \xi(t)=0\}.
$$
Then we  set $i_t=\i_{\lfloor \hat\zeta(t+1)X_t \rfloor}$ if $X_t\ge 1$, and $i_t=\i_{\lfloor \hat\zeta(t+1)n \rfloor}$ otherwise. Now,
$$
\xi_j(t+1)=\begin{cases}
\xi_j(t), &\text{if } j\notin\{i_t-1,i_t,i_t+1\};\\
\zeta_{j-i_t}(t+1), &\text{if }  j\in\{i_t-1,i_t,i_t+1\}.
\end{cases}
$$
We  define the sigma-algebra
$\F_t=\sigma\left(\bar\zeta(1),\bar\zeta(2),\dots,\bar\zeta(t)\right)$, then $\xi(t)$ is $\F_t$-measurable.

Since the Markov chain $\xi(t)$ is irreducible, aperiodic  and on a finite state space, it converges to the unique stationary distribution $\pi^{(n)}$. Let 
$$
\nu^{(n)}(p)=\sum_{\xi\in \{0,1\}^n}  \pi^{(n)}(\xi)\, {\bf 1}_{\xi_i=1}
$$ 
be the probability that vertex $i$ has value $1$ in this stationary distribution; by symmetry, this quantity does not depend on $i$ and equals the expected number of ones under the stationary distribution, divided by $n$.

Now let $n\to\infty$. It is not hard to guess {\it intuitively} that $\lim_{n\to\infty}\nu^{(n)}(p)=1$ if $p>2/3$, since every time we replace at least $1$ zero with on average $3(1-p)$ zeros. We shall formulate the following 
\begin{conj}
\begin{itemize}
\item[(a)] There exists $\nu(p)=\lim_{n\to\infty}\nu^{(n)}(p)\in[0,1]$.
\item[(b)] The function $\nu(p)$ is monotone increasing in $p\in[0,1]$.
\item[(c)] There exists $p_c\in [0,1]$ such that $\nu(p)=1$ for $p>p_c$ and $\nu(p)<1$ for $p<p_c$.
\item[(d)] $p_c$ lies strictly between $0$ and $1$.
\end{itemize}
\end{conj}
None of the parts of the above conjecture are rigorously shown. At the same time~\cite[Theorem~2.1]{BK} claims that $\lim_{n\to\infty} \nu^{(n)}(p)=1$ for $p>p_*=0.54\dots$. Hence, if $p_c$ exists, then  $p_c<p_*$. At the same time it appears much harder to show that $p_c>0$, or at least that $\limsup_{n\to\infty} \nu^{(n)}(p)<1$ for small enough $p$. The latter statement was eventually proven in~\cite{MZ}, using a very involved combinatorial {\em avalanche} method.

While the statement of \cite[Theorem~2.1]{BK} is correct, its proof is, unfortunately, not rigorous. The authors correctly compute the probabilities of bounds of possible changes in $D_t$, which the size of the smallest contiguous area containing all zeros (precisely defined later in our article, see~\eqref{Ddef}), and show that  its drift is negative whenever $p>p_*$.  Then they proceed to compare $D_t$ with  Markov chain $X_t$ with some ``holding'' probabilities $\alpha_i$, and claim that $\lim_t\E X_t \ge \lim_t \E D_t$ as long as ``alphas are defined  to maximize  $\E X_t$ subject to some relevant constraints'' (just above Lemma~2.1 in~\cite{BK}). This cannot in general be true, as, for example, $\lim_{t\to\infty}\E D_t$ might not even exist. It is also not clear to us how exactly the authors derived Lemma~2.2 from~\cite{FMM}.

The purpose of this short note is twofold. Firstly, we show how one can make the results of~\cite{BK} rigorous. Secondly, we obtain a better bound on the critical probability by showing that $\nu(p)=1$ for all $p>p_\m$ where $p_\m=0.45\dots<p_*$. This is done in the next section. The final section contains the required statements for stochastic processes with drift, which are crucial in order to make the proof mathematically rigorous.

For some other recent results on Bak-Sneppen model please see~\cite{BB,BR,FRA,GKW,MS} and references therein.

\section{Main result}\label{main}
This is the main result of our paper.
\begin{thm}\label{t1}
Let $p_\m=0.45\dots$ be the only real positive solution of the equation
\begin{align}\label{eqsol}
p^5+4p^4+2p^3+3p^2=1.
\end{align}
Then for all $p\in(p_\m,1]$ there exists $\nu(p)=\lim_{n\to\infty} \nu^{(n)}(p)=1$. Hence, if $p_c$ exists, $p_c\le p_\m$.
\end{thm}
\begin{rema}
Note that the simulations suggest that $p_c\approx 0.36$
(see, e.g.~\cite[Figure 1]{BK}).
\end{rema}
Similarly to~\cite{BK}, we will define $X_t=|\{i:\ \xi_i(t)=0\}|$ as the total number of zeros at time $t$. The proof of Theorem~\ref{t1} will be based on finding some sort of a Lyapunov function $M_t$, which is a function of $\xi(t)$, and is thus $\F_t$-measurable, satisfying the following properties:
\begin{itemize}
\item $M_t\ge X_t\ge 0$ for all realizations;
\item $M_t$ has uniformly bounded up-jumps;
\item  there are some constants $C>0$ and $\eps>0$ such that  on the event $\{M_t\ge C\}$ we have $\E(M_{t+1}-M_t \| \F_t,M_{t+1}\ne M_t)\le -\eps $;
\item the probability  $\P(M_{t+1}\ne M_t\|\F_t)$ is bounded below by $\frac{\rm const}{M_t+1}$.
\end{itemize}
Once such $M_t$ is found (see Lemma~\ref{lem_drift}), Lemma~\ref{lemrenew} (which itself follows from Lemma~\ref{lemFMM}) would imply that, loosely speaking,  $\E M_t$ remains ``on average'' bounded by a quantity independent of $n$, and hence the same holds for $X_t$. Since $\E X_t$ in fact converges as $t\to\infty$, we will get the desired result.
\begin{rema}
The same proof can be used almost verbatim to show \cite[Theorem 2.1]{BK} rigorously; the only difference is in the computation of the expected drift of $D_t$ (defined below), instead of~$M_t$.
\end{rema}

As in~\cite{BK}, let $D_t$ be the diameter of the smallest configuration containing all the zeros. Formally, for a given configuration $\xi(t)=(x_0,\dots,x_{n-1})$ let
$$
{\bf Z}_t=(i,i+1,\dots,j-1,j)
$$
be  a sequential subset of indices of $(0,1,2,\dots,n-1)$ with periodic boundary conditions, satisfying the following two properties
\begin{itemize}
\item[(a)] $\xi_k(t)=1$ for all $k\not\in {\bf Z}_t$;
\item[(b)] ${\bf Z}_t$ has the smallest number of elements amongst all such subsets.
\end{itemize}
If there is more than one such subset at time $0$, choose any of them arbitrarily\footnote{Example: $\xi(t)=(1,0,0,1,0,0)$, then ${\bf Z}_t$ is either $(1,2,3,4,5)$ or $(4,5,0,1,2)$.}.
Now we can define
\begin{align}\label{Ddef}
D_t&=\card({\bf Z}_t) \ \ge X_t.
\end{align}
Also let $(l_t,r_t)$ be the pair  with the first and the last index of ${\bf Z}_t$, and set $l_t=r_t=0$ if $X_t=n$ or $X_t=0$, i.e.\ there are no ones or no zeros amongst $\xi_i(t)$. Note that if $\xi_t=(x_0,\dots,x_{n-1})$ and $1\le X_t\le n-1$ then
$$
x_{l_t-1}=1,\quad x_{l_t}=0,\quad x_{r_t}=0,\quad x_{r_t+1}=1.
$$

Suppose that $D_t\ge 6$, and for some positive constant 
$
\b\in(0,1/2)
$
define a ``corrected'' diameter of the configuration
\begin{align}\label{eqM}
M_t=M(\xi(t))=D_t+1 - \b\left ({\bf 1}_{x_{l_t+1}=0}+{\bf 1}_{x_{r_t-1}=0}\right)\in (D_t,D_t+1]
\end{align}
that is, $M_t$ differs from $D_t$ by at most $1$. For definiteness, let $M_t=0$ whenever $D_t<6$.

Note that if one of the ``deeply'' internal indices\footnote{namely, $l_t+2< i_t< r_t-2$} of ${\bf Z}_t$ is chosen, then $x_{l_t},x_{l_t+1},x_{r_t-1},x_{r_t}$ do not change  and  ${\bf Z}_{t+1}={\bf Z}_t$, e.g.\ 
$$
\dots {\bf 111}\ 00\blue{1\underline{0}1}10\ {\bf 111}\dots 
\longrightarrow
\dots {\bf 111}\ 00\blue{***}10\ {\bf 111}\dots 
\quad
({\bf Z}_{t+1}= {\bf Z}_t,\ l_{t+1}=l_t,\ r_{t+1}=r_t).
$$
On the other hand, if $i_t\in\{l_t,l_t+1,r_t,r_t-1\}$ then ${\bf Z}_t$ might change; it can increase by at most one point; if $i_t\in\{l_t+2,r_t-2\}$ then ${\bf Z}_t$ does not change but $M_t$ still can.

The above analysis fails, however, if one can make $\card({\bf Z}_{t+1})$ even smaller by flipping indices, which can happen when one of the  ``deeply'' internal zeros is chosen, e.g.
\begin{align*}
\xi(t)=(1,0,1,\blue{0,\underline{0},0},1,0,1)
&\longrightarrow
\xi(t+1)=(1,0,1,\blue{1,1,1},1,0,1)\\
{\bf Z}_t=\{1,2,3,4,5,6,7\}&\longrightarrow {\bf Z}_{t+1}=\{7,8,0,1\}.
\end{align*}
In this case $l_{t+1}=r_t$, $r_{t+1}=l_t$ and $D_{t+1}\le D_t-1$ and thus $M_{t+1}\le M_t$.
Hence, when $t\ge 1$, in case of ties for the choice of ${\bf Z}$, we will be able to choose $l_t$ and $r_t$ such that at least one of the following events hold:
\begin{align}\label{eqtie}
l_t=l_{t-1}\text{ or }r_t=r_{t-1}\text{ or }
\{l_{t+1}=r_t, r_{t+1}=l_t\}.
\end{align}

\begin{defn}
We say that the indices $l$ and $r$ {\em flip} when the configuration changes from $\xi(t)$ to $\xi(t+1)$ such that the last event in~\eqref{eqtie} holds\footnote{Example: $n=7$, $\xi(t)=(1,0,0,0,0,0,1)$, $(l_t,r_t)=(1,5)$, $D_t=5$. Now suppose the middle zero is chosen, and, together with its both neighbours, replaced by ones. The new configuration shall be
$\xi(t+1)=(1,0,1,1,1,0,1)$, $(l_{t+1},r_{t+1})=(5,1)$, $D_{t+1}=4$.}. 
\end{defn}
Note that if $l$ and $r$ do not flip, then 
\begin{align*}
{\bf Z}_{t+1}\subseteq {\bf Z}_t\quad\text{or}\quad
{\bf Z}_{t+1}={\bf Z}_{t} \cup \{l_t-1\}\quad\text{or}\quad 
{\bf Z}_{t+1}={\bf Z}_{t} \cup \{r_t+1\}
\end{align*}
according to the above arguments.

\begin{lemma}\label{lem_flip}
Suppose that $D_t\ge 6$, and suppose that the indices $l$ and $r$ flip between times $t$ and $t+1$. Then
$$
D_{t+1}\le D_t-1, \quad\Longrightarrow \quad M_{t+1}\le M_t-(1-2\b)<M_t.
$$
\end{lemma}
\begin{proof}
The first statement follows from the fact that the flip occurs only if $\card({\bf Z}_t)$ decreases; the second from the fact that $2\b<1$.
\end{proof}

Let $l=l_t$, $r=r_t$. When $D_t\ge 6$ define
$$
{\bf E}_t=\{l,l+1,l+2\}\cup\{r,r-1,r-2\}\subseteq {\bf Z}_t
$$
as the set indices of the three+three points at both ends of the zero area.

Recall that $i_t$ denotes the index of the zero, chosen to be replaced with both its neighbours, and note that ${\bf 1}_{i_t\in{\bf E_t}}$ is a function of $\xi(t)$ and $\hat\zeta(t+1)$.

\begin{lemma}\label{lem_inside}
Suppose that $D_t\ge 6$. Then $M_{t+1}\le M_t$ on $i_t\not\in {\bf E}_t$.
\end{lemma}
\begin{proof}
On the event described in the statement, either $D_{t+1}=D_t$ (and there were no changes of the configuration near its endpoints, hence $M_{t+1}=M_t$), or the indices $l$ and $r$ flipped, and in this case the result follows from Lemma~\ref{lem_flip}.
\end{proof}

\begin{lemma}\label{lem_drift}
Suppose that $p>p_\m$ and $\beta$ is given by~\eqref{eqbeta}. Then \begin{align*}
\Delta_{t+1}:=\E(M_{t+1}-M_t\| \F_t)&\le 0
\quad\text{on the event }\{M_t\ge 8\}.
\end{align*}
Moreover, there is an $\eps=\eps(p)>0$, depending on $p$ only, such that 
\begin{align*}
\Delta_{t+1}&\le -\eps
\quad\text{on the event }\{M_t\ge 8\} \cap
\{i_t\in {\bf E}_t\}.
\end{align*}
\end{lemma}
\begin{rema}
The inequalities above are exactly of the type of those established in~\cite{BK} for the quantity $D_t$. By slightly modifying $D_t$ with extra terms in~\eqref{eqM}, we will be able to obtain a  better bound on $p_c$.
\end{rema}

\begin{proof}[Proof of Lemma~\ref{lem_drift}]
Throughout the proof we write $\xi(t)=x=(x_0,x_1,\dots,x_{n-1})$, $l=l_t$ and $r=r_t$.
Also, because of Lemma~\ref{lem_flip}  for the rest of the proof we may assume that the flip of the indices $l$ and $r$ does not happen (should that occur, $M_t$ decreases by at least $1-2\b>0$).
Observe also that $M_t\ge 8$ ensures that $D_t\ge 6$, and hence~\eqref{eqM} holds. 

For the moment,  assume that $D_t\le n-2$, the remaining two cases will be investigated later. Due to the symmetry, it is sufficient to study only the left end of the zero configuration; the drift on the right end is identical. Also w.l.o.g.\ assume $l=2$ (and hence $x_0=x_1=1$, $x_2=0$).

There are four possibilities for the beginning of the configuration $x=\xi(t)$:
\begin{multicols}{2}
\begin{itemize}
\item[(a)] $x=(1,1,0,{\bf 0,0},\dots)$;
\item[(b)] $x=(1,1,0,{\bf 0,1},\dots)$;
\item[(c)] $x=(1,1,0,{\bf 1,0},\dots)$;
\item[(d)] $x=(1,1,0,{\bf 1,1},\dots)$.
\end{itemize}
\end{multicols}

The following calculations are done by thoroughly examining the 8 possible cases where ``*0*'' is replaced by ``000'',``001'', $\dots$, ``111'' respectively, and the probability of choosing ``1'' is~$p$.

Conditioned on choosing one of the three zeros shown in case (a), each with equal probability, the drift $\Delta_{t+1}$ is bounded above by
$$
T_{00}=-\frac13\left[p^3+p^2+p-1+\b(1-p)\right]
-\frac13\left[ p(p^2+p+1)+\b(p+1)(1-p)^2 \right]
-\frac{\b(1-p)}3  +\beta.
$$
In case (d), conditioned on choosing the shown zero, $\Delta_{t+1}$ is  bounded above by
$$
T_{11}=-(2p-1)(p^2+p+1)-\beta(1-p)^2(1+p).
$$
In case (b), conditioned on choosing one of the two zeros, $\Delta_{t+1}$ is bounded above by
$$
T_{01}=\frac{T_{11}}{2}-\frac12\left[p(p^2+p+1)+\beta(1+p)(1-p)^2 \right]
+\b.
$$
Finally, in case (c), conditioned on choosing one of the two zeros, $\Delta_{t+1}$ is  bounded above by
\begin{align}
\label{eqT10}
T_{10}=-\frac 12\left[p^3+p^2+p-1+\beta(1-p)\right]-\frac{\beta(1-p)}2.
\end{align}
The reason we have inequalities  for the drift, rather than equalities, is that sometimes we do not know whether the new configuration ${\bf Z}_{t+1}$ starts with ``00'' or ``01'', and we take the worst case scenario (i.e.\ ``01'').

While the computations of $T_{ij}$ are quite tedious, for the sake of completeness we will present the detailed calculation of $T_{10}$ (case (c)). Suppose w.l.o.g that $l_t=2$, and thus $\xi(t)=11010??\dots$ where question marks correspond to the unknown values. Then $M_t=R_t-L_t+1$ where $L_t=l_t$ and $R_t=r_t-\b {\bf 1}_{x_{r_t-1}=0}$. Also, $M_{t+1}=R_t-L_{t+1}+1$ if $i_t=2$ or $i_t=4$, and we need to compute the value of $L_{t+1}$. In the table below we present the values of $L_t-L_{t+1}$ depending on the values of $\bar\zeta_t$ as well as the probabilities of those outcomes.
\\[5mm]
If $i_t=2$, then:
\begin{tabular}{l|l|l}
Configuration & $L_t-L_{t+1}$ & Probability\\
\hline
$1{\bf 00*}0??\dots$ & $1-\b$ & $(1-p)^2$ \\
$1{\bf01*}0??\dots$ & $1$ & $p(1-p)$ \\
$1{\bf100}0??\dots$ & $-\b$ & $p(1-p)^2$ \\
$1{\bf101}0??\dots$ & $0$ & $p^2(1-p)$ \\
$1{\bf110}0??\dots$ & $-1-\b$ & $p^2(1-p)$ \\
$1{\bf111}0??\dots$ & $-2$ & $p^3$\\ 
\hline
\end{tabular}
\\[5mm]
If $i_t=4$, then:
\begin{tabular}{l|l|l}
Configuration & $L_t-L_{t+1}$ & Probability\\
\hline
$110{\bf0**}?\dots$ & $-\b$ & $1-p$ \\
$110{\bf1**}?\dots$ & $0$ & $p$\\
\hline
\end{tabular}
\\[5mm]
Here $*$ stands for the unimportant value which is generated in $\xi(t+1)$.
Since $i_t=2$ or $4$ with equal probibility, we get
$$
T_{10}=\frac{1}{2}\left[(1-\b)(1-p)^2+p(1-p)-\b p(1-p)^2+(-1-\b)p^2(1-p)-2p^3\right]
+
\frac{1}{2}\left[(-\b)(1-p)\right]
$$
which is the same as~\eqref{eqT10}.
\\[5mm]
Next, we need to ensure that  $T(p,\b):=\max_{i,j\in\{0,1\}} T_{ij}$ is strictly negative. It can be shown using elementary algebra that 
$$
\frac{\partial T_{ij}}{\partial p}<0\text{ for all }i,j=0,1,\quad p\in(0,1),\quad \b\in(0,1/2).
$$
Hence, if $T(\tilde p,\b)=0$ for some $\tilde p$, then
$T( p,\b)< 0$ for all $p> \tilde p$.

By examining all the cases, we find that the largest value of $p$ for which $T(p,\b)$ can be made non-positive, is the real solution of the equation~\eqref{eqsol}; in this case
\begin{align}\label{eqbeta}
\b&=\b_\m=-\frac{4p_\m^3+p_\m^2+p_\m-2}{2p_\m^3-2p_\m^2+3}=0.34656\dots,
\\
T_{00}=T_{11}&=0,\quad T_{01}=-0.0669\dots,\quad T_{10}=-0.106\dots.
\nonumber
\end{align}
Hence, for every $p>p_\m$ we have $T(p,\b_\m)<-\eps(p)<0$.

Finally, the cases where $D_t=n$ or $D_t=n-1$ are trivial, 
since  $\E (M_{t+1}-M_t\|\F_t)=\E (D_{t+1}-D_t\|\F_t)<-[1-(1-p)^3]<0$ on~${\bf E}_t$  in the first case $x=(0,0,\dots,0)$, and
\begin{align*}
\E (M_{t+1}-M_t\|\F_t)&=\frac 12\left[  (+1)p+(-1)2p^2(1-p)+(-2)p^2\right]
\\ & +
\frac 12\left[(-1)p(1-p)^2+(-2)p(1-p)+(-3)p  \right]=-\frac{p(5-p^2)}{2}<0
\end{align*}
on~${\bf E}_t$ in the second case where e.g.\ $x=(0,0,1,0,0,\dots,0)$.
\end{proof}
\vskip 5mm
Now we are ready to present the proof of the main result.
\begin{proof}[Proof of Theorem~\ref{t1}]
Let $M_t$ be defined by~\eqref{eqM}, and let the sequence of stopping times $\tau_k$ be defined by $\tau_0=0$ and
$$
\tau_{k+1}=\inf\{t>\tau_k:\ i_t\in {\bf E}_t\}.
$$
Then $M_t$ and $\tau_k$ satisfy the conditions of Lemma~\ref{lemrenew} with $\tilde\F_t=\sigma\left(\F_t,{\bf 1}_{i_t\in{\bf E}_t}\right)$, $C=8$, $b=2$, $r=1$ and $\eps$ given by Lemma~\ref{lem_drift}. Indeed, 
$M_{t+1}\ne M_t$ if and only if either there is a flip of $l$ and $r$ (in which case $M_{t+1}<M_t$ by Lemma~\ref{lem_flip}), or if $i_t\in {\bf E}_t$. Since each zero in the configuration~$\xi(t)$ is chosen with equal probability, the probability of the latter event is bounded below as follows:
$$
\P(i_t\in {\bf E}_t\|\F_t)= \frac{\card\{j\in {\bf E}_t:\ \xi_j(t)=0 \}}{X_t}\ge \frac{1}{D_t} \ge \frac{1}{M_t+1}
$$
from~\eqref{eqM} since $\b<1/2$. Hence $\tau_{k+1}-\tau_k$ is bounded above by a geometric random variable with expectation $M_{\tau_k}+1$.

Now, Lemma~\ref{lemrenew} implies that
\begin{align}\label{eqCes}
\limsup_{T\to\infty}\frac{  \sum_{t=1}^T \E M_t}{T}\le \tilde R
\end{align}
for some $\tilde R$ not dependent on $n$. At the same time $M_t$ is a function of a positive recurrent finite-state Markov chain, hence  as $t\to\infty$
$$
\E M_t\to \sum_{\xi\in \{0,1\}^n}  \pi^{(n)}(\xi)\, M(\xi)=\tilde \E M(\xi)=:\mu^{(n)}
$$ 
where $\tilde\E$ denotes the expectation under the stationary measure $\pi^{(n)}$ for the chain $\xi$, and  $\mu^{(n)}$ is the expectation of $M$ under this measure. Therefore, by Ces\`aro summation from~\eqref{eqCes} we obtain that $\mu^{(n)}\le \tilde R$. Since $X_t\le D_t\le M_t+1$, the expected fraction of zeros in the limit is given by
\begin{align}\label{eqfin}
1-\nu^{(n)}(p)=\frac{\tilde\E X_t}{n}
\le \frac{\tilde\E D_t}{n} \le \frac{\tilde\E (M_t+1)}{n} 
\le \frac{\tilde R +1}{n},
\end{align}
The RHS of~\eqref{eqfin} converges to zero as $n\to\infty$, yielding the statement of the theorem.
\end{proof}

\begin{rema}
By analyzing a larger (but still a finite) set of cases, one can even get a better upper estimate for $p_c$, namely $p_c\le 0.419533$, which is closer to the estimated value of $0.36$. In order to do this, one can introduce the following more subtle supermartingale (on $D_t\ge 8$)
\begin{align*}
\tilde M_t=
D_t
&- \a\left ({\bf 1}_{(x_{l+1},x_{l+2})=(0,0)}+{\bf 1}_{(x_{r-1},x_{r-2})=(0,0)}\right)
\\ &
- \b\left ({\bf 1}_{(x_{l+1},x_{l+2})=(1,0)}+{\bf 1}_{(x_{r-1},x_{r-2})=(0,1)}\right)
\\ &
- \g\left ({\bf 1}_{(x_{l+1},x_{l+2})=(0,1)}+{\bf 1}_{(x_{r-1},x_{r-2})=(1,0)}\right)
\end{align*}
where $\a=0.3764287$, $\b=0.078811$, and $\g=0.423494$. For these particular values, we get that $\eps(p)$ used in Lemma~\ref{lem_drift} is at least $3.6\times 10^{-8}$.
\end{rema}

\section{Appendix}
\begin{lemma}\label{lemFMM}
Consider a real-valued non-negative process $Y_t$, $t=0,1,2,\dots$ adapted to the filtration~$\F_t$ for which there exist $C,b,\eps>0$ such that the differences $\Delta_{t+1}=Y_{t+1}-Y_t$ satisfy
\begin{itemize}
\item[(a)] $\Delta_{t+1}\le b\text{ for all }t;$
\item[(b)] $\E\left(\Delta_{t+1} \cdot {\bf 1}_{\Delta_{t+1}\ge -b}\|\F_t\right)\le -\eps\text{ on } Y_t\ge C.$
\end{itemize}
Then there exists an $R_m=R_m(C,b,\eps)<\infty$ such that for every $m=1,2,\dots$
\begin{align*}
\limsup_{t\to\infty} \E \left(Y_t^m\right)\le R_m.
\end{align*}
\end{lemma}
\begin{proof}
The idea is borrowed from the proof of~\cite[Theorem 2.1.7]{FMM}.
Let $0<h<1/b$. Since $e^x\le 1+x+x^2$ for $|x|\le 1$,  we have
\begin{align}\label{eqexp}
e^{h\Delta_t}&\le \begin{cases}
1+ h\Delta_t+ h^2\Delta_t^2, &\text{if } \Delta_t \ge -b;
\\
e^{-hb},&\text{otherwise}
\end{cases} \nonumber \\
& \le  \left(1 -e^{-hb}+ h^2 b^2 +h\Delta_t\right) {\bf 1}_{\Delta_t\ge -b}+
e^{-hb}
\end{align}
since $\Delta_t^2\le b^2$ on $\{\Delta_t\ge -b\}$.
Then, taking the expectation of~\eqref{eqexp}, we have on $\{Y_{t-1}\ge C\}$
\begin{align*}
\E\left(e^{h\Delta_t}\|\F_{t-1}\right)
&\le  
\left(1 -e^{-hb}+ h^2 b^2 \right) \P(\Delta_t\ge -b)+
e^{-hb}-h\eps\\
&\le 
1 -e^{-hb}+ h^2 b^2 +e^{-hb}-h\eps=1+ h(h b^2 -\eps)
\le 1-h\eps/2=:\g<1
\end{align*}
by taking  $h$  so small that $hb^2<\eps/2$. 
Consequently,
\begin{align}\label{eqFMM}
\E\left(e^{hY_t}\right)&=
\E\left(e^{hY_{t-1}}\E\left(e^{h\Delta_t}\|\F_{t-1}\right)\right)
\le 
\E\left(e^{hY_{t-1}}  \left[ 
\g {\bf 1}_{Y_{t-1}\ge C}+ e^{hb} {\bf 1}_{Y_{t-1}< C}\right] 
\right)\nonumber \\
&\le \g \E\left(e^{hY_{t-1}} \right)+ e^{h(C+b)}.
\end{align}
Now, by re-iterating~\eqref{eqFMM}, we obtain that
$$
\E\left(e^{hY_t}\right)\le e^{h(C+b)} (1+\g+\dots+\g^{t-1})+\g^t \E\left(e^{hY_0}\right)<
\frac{e^{h(C+b)}}{1-\g}+\g^t \E\left(e^{hY_0}\right).
$$
Hence, since $e^{hy}\ge (hy)^m/m!$ for all $y,h\ge 0$, and $m=1,2,\dots$,
$$
\limsup_{t\to \infty} \E \left(Y_t^m\right)\le \limsup_{t\to \infty} \frac{m! \, \E\left(e^{hY_t}\right)}{h^m}
\le \frac{m!\,e^{h(C+b)}}{h^m(1-\g)}=\frac{2\, m!\,e^{h(C+b)}}{\eps\, h^{m+1} }<\infty.
$$
\end{proof}

The proof of the next statement uses some ideas from the renewal theory.
\begin{lemma}\label{lemrenew}
Consider a real-valued non-negative process $M_t$, $t=0,1,2,\dots$ adapted to the filtration~$\tilde\F_t$ for which there exist $C,b,r,\eps>0$ such that
\begin{enumerate}
\item[(a)] $M_{t+1}-M_t\le b\text{ for all }t;$
\item[(b)] $\E\left[\left(M_{\tau_{k}+1}-M_{\tau_{k}}\right)\cdot{\bf 1}_{M_{\tau_{k}+1}-M_{\tau_{k}}\ge-b}\|\tilde\F_{\tau_k}\right] \le -\eps\text{ on }M_{\tau_k}\ge C;$
\item[(c)] $M_{t+1}\le M_t \text{ on } t\not\in\cup_k\tau_k$;
\item[(d)]  $\E( \tau_{k+1}-\tau_k\|\tilde\F_{\tau_k})\le r (1+M_{\tau_k})$
\end{enumerate}
for some strictly increasing positive integer sequence of stopping times $\tau_k$, $k=0,1,2\dots$.
Then
\begin{align*}
\limsup_{T\to\infty}  \frac{\sum_{t=1}^T \E M_t}T\le \tilde R
\end{align*}
where $ \tilde  R$ depends on $C,b,r,\eps$ only.
\end{lemma}
\begin{proof}
First, let $Y_k=M_{\tau_k}$. Because of (a), (b) and (c), $Y_k$ satisfies the conditions of Lemma~\ref{lemFMM}. 

Let $N(t)=\max\{k:\ \tau_k\le t\}$ be the counting process. Then, from (c),
\begin{align}\label{eqsum}
\sum_{t=1}^T M_t\le \sum_{t=1}^{\tau_{N(T)+1}	} M_t
\le
\sum_{k=0}^{N(T)+1} M_{\tau_k} (\tau_{k+1}-\tau_k)
\le \sum_{k=1}^{T+1} M_{\tau_k} (\tau_{k+1}-\tau_k)
\end{align}
since $N(t)\le t$. On the other hand, from (d),
\begin{align*}
\E\left(M_{\tau_k} (\tau_{k+1}-\tau_k)\right)&=
\E\left[M_{\tau_k} \E(\tau_{k+1}-\tau_k\|\tilde\F_{\tau_k})\right]
\le \E\left[r M_{\tau_k} (1+M_{\tau_k})\right]
= r\,\E\left[ Y_k +Y_k^2\right]
\end{align*}
so that
 $$
 \limsup_{k\to\infty}\E\left(M_{\tau_k} (\tau_{k+1}-\tau_k)\right)
 \le r (R_1+R_2)
 $$
where $R_1,R_2<\infty$ depend only on $C,b,\eps$ by Lemma~\ref{lemFMM}. 
This, together with~\eqref{eqsum}, implies  that
$$
\limsup_{T\to\infty}\frac{\displaystyle  \sum_{t=1}^T \E M_t}{T}\le 
\limsup_{T\to\infty}\frac{\displaystyle \sum_{k=0}^{T+1} \E\left( M_{\tau_k} (\tau_{k+1}-\tau_k)\right)}{T+2}\cdot \frac{T+2}{T}
\le  r (R_1+R_2)
$$
as the Ces\`aro mean.
\end{proof}


\end{document}